\documentclass[reqno,10pt]{amsart}
\usepackage{amscd,amssymb,amsmath,color,url}
\usepackage{latexsym}
\usepackage[all]{xy}
\usepackage{defs}
\usepackage[colorlinks=true]{hyperref}

\def\..{{,\dots ,}}

\def\Br{{\rm Br}}
\def\Ram{{\rm Ram}}

\def\rmlog{{\rm log}}

\def\di{{\diamond}}

\def\slope{{\rm slope}}
\def\hyp{{\rm hyp}}

\begin{document}

\title{Wild coverings of Berkovich curves}

\author{Michael Temkin}
\address{Einstein Institute of Mathematics, The Hebrew University of Jerusalem, Giv'at Ram, Jerusalem, 91904, Israel}
\email{temkin@math.huji.ac.il}

\maketitle

\section{Introduction}
This paper is an extended version of the author's talk given at the conference "Non-Archimedean analytic geometry: theory and practice" held in August 2015 at Papeete, and I wish to thank the organizers. It gives a brief overview of results and methods of works \cite{CTT} and \cite{radialization} on the structure of finite morphisms between Berkovich curves.

The structure of tame morphisms between smooth Berkovich curves is pretty well-known and it is completely controlled by the simultaneous semistable reduction theorem, see, for example, \cite{ABBR}. The structure of wild morphisms was for a long time terra incognita, though one should mention some special results recently obtained by Faber in \cite{XFaber1} and \cite{XFaber2}. In this project we obtain a relatively complete description of the combinatorial structure of an arbitrary finite morphism $f\:Y\to X$ between smooth Berkovich curves. It is divided into two parts.

\subsection{The different function}\label{diffunsec}
In a joint work \cite{CTT} with A. Cohen and D. Trushin we study the {\em different function} $\delta_f\:Y^\hyp\to[0,1]$ that assigns the different $\delta_{\calH(y)/\calH(f(y))}$ to a point $y$ of type 2, 3 or 4. In other words we study the analytic behavior of the most important invariant that measures wildness of an extension of valued fields, the different. It turns out that $\delta_f$ controls the minimal semistable model of $f$, and a balancing condition for the slopes of $\delta_f$ at a type 2 point $y\in Y$ extends the local Riemann-Hurwitz formula to the wild case.

\subsection{The multiplicity function and radialization}\label{multsec}
Both $X$ and $Y$ have canonical exponential metrics and $f$ is piecewise monomial with respect to them. The behavior of $f$ as a piecewise monomial function is controlled by the {\em multiplicity function} $n_f$ that assigns $[\calH(y):\calH(f(y))]$ to $y$. This function and the multiplicity loci $N_{f,\ge d}=\{y\in Y\vert\ n_f(y)\ge d\}$ are described in \cite{radialization}. In particular, it is shown that $n_f$ is radial with respect to a large enough skeleton $\Gamma_Y\subset Y$. A central player in this study is a {\em profile function} $\phi_f\:Y^\hyp\to P_{[0,1]}$ encoding the radii of all sets $N_{f,\ge d}$ around the skeleton, where $P_{[0,1]}$ is the set of piecewise monomial bijections of $[0,1]$ onto itself. Furthermore, $\delta_f$ can be retrieved from $\phi_f$ via composing with a character $P_{[0,1]}\to[0,1]$ and $\phi$ is an analytic family of the classical Herbrand functions.

\subsection*{Acknowledgments}
I would like to thank the referee for pointing out various inaccuracies in the first version of the paper.

\section{Semistable reduction and tame morphisms}
In this section we summarize the relatively well-known properties of curves and morphisms.

\subsection{Conventions}

\subsubsection{Ground field}
We work over a fixed algebraically closed non-archimedean {\em analytic} (i.e. real-valued complete) field $k$ with a non-trivial valuation.

\subsubsection{Nice curves}
By a {\em nice curve} we mean a rig-smooth connected separated compact $k$-analytic curve.

\subsubsection{Subgraphs}
By a {\em subgraph} $\Gamma$ of a nice curve $C$ we mean a connected topological subgraph $\Gamma\subset C$ with finitely many vertices and edges such that the set of vertices $\Gamma^0$ consists of points of $C$ of types 1 and 2 and contains at least one point of type 2.

\subsection{Semistable reduction for curves}

\subsubsection{Skeletons of curves}
A subgraph $\Gamma$ is called a {\em skeleton} of $C$ if $C\setminus\Gamma^0$ is a disjoint union of open discs $D_i$ and semi-annuli $A_1\..A_n$ (i.e. $A_i$ is either an open annulus or a punched open disc) and the edges of $C$ are the skeletons of $A_1\..A_n$. The following skeletal version of the semistable reduction theorem is easily seen to be equivalent to its classical versions.

\begin{theor}
Any nice curve possesses a skeleton.
\end{theor}

\subsubsection{Combinatorial structure of the curve}
In a sense, a skeleton of a curve provides the best possible combinatorial description of the curve. In particular, the complement $C\setminus\Gamma$ of a skeleton is a disjoint union of discs and there is a canonical deformational retraction $q_\Gamma\:C\to\Gamma$.

\subsubsection{Genus}
For any point $x\in C$ we define the {\em genus} $g(x)$ to be the genus of $\wHx/\tilk$ if $x$ is of type 2 and zero otherwise. The {\em genus} of $C$ is then defined to be $g(C)=h^1(C)+\sum_{x\in C}g(x)$; it is finite and equals to $g(C)=h^1(\Gamma)+\sum_{v\in\Gamma^0}g(v)$. This gives the usual genus of an algebraic curve when $C$ is proper, but $g(C)$ is a meaningful invariant for nice curves with boundary too.

\subsubsection{Exponential metric}
Let $A$ be an open or closed annulus, i.e. $A$ is isomorphic to the domain in $\bfA^1_k$ given by $r_2<|t|<r_1$ or $r_2\le|t|\le r_1$. The number $r(A)=r_1/r_2$ depends only on $A$ and it is called the {\em radius} of $A$. Given an interval $I\subset C$ one defines its {\em radius} (or exponential length) by $r(I)=\sup \prod_{i=1}^nr(A_i)$, where the supremum is taken over all finite sets of disjoint open annuli $A_i\subseteq C$ such that the skeleton of each $A_i$ lies in $I$. It turns out that $r$ defines an exponential metric on $C$ whose singular points are precisely the points of type 1. In other words, $r([a,c])=r([a,b])r([b,c])$ for an interval $[a,c]\subset C$ with a point $b\in[a,c]$, and $r([a,b])=\infty$ if and only if the set $\{a,b\}$ contains a point of type 1.

\begin{rem}
(i) If $I$ is the skeleton of $A$ then $r(I)=r(A)$. In fact, this is the main property of the radius one should check in order to establish all other properties.

(ii) We prefer to work with the exponential metric in this paper, but one often considers its logarithm, which is a usual metric. For example, if $I$ is the skeleton of an annulus $A$ then the length of $I$ is the modulus of $A$. The classical metric is only canonical up to rescaling since its definition involves a choice of the base of the logarithm.
\end{rem}

\subsubsection{Radius parametrization}
If $I$ is closed with an endpoint $a$ of type different from 1 then $x\mapsto r([a,x])^{-1}$ provides the canonical homeomorphism $I=[r(I)^{-1},1]$ that we call {\em radius parametrization} of $I$. Note that $I=[0,1]$ if and only if the second endpoint is of type 1.

\begin{rem}
(i) If $E$ is a unit disc with maximal point $q$ then for any point $x\in E$ there exists a unique interval $[x,q]$ and $r(x)=r([x,q])^{-1}$ is the usual radius function of the disc.

(ii) In the same way, any skeleton $\Gamma$ induces a radius function $r_\Gamma\:C\to[0,1]$ that measures the inverse exponential distance to the skeleton.
\end{rem}

\subsubsection{Enhanced skeleton}
We naturally enhance a skeleton $\Gamma$ of a curve to a {\em metric genus graph} in which each vertex is provided with a genus and each edge is provided with a radius
(exponential length).

\subsection{Semistable reduction for morphisms}

\subsubsection{Morphisms and metrics}
Let $f\:Y\to X$ be a non-constant morphism of nice curves. It is easy to see that $f$ is {\em pm} or piecewise monomial in the sense that for each interval $I\subset Y$ the set $f(I)$ is a graph and the map $I\to f(I)$ is pm with integral slopes with respect to the radii parameterizations. Moreover, the multiplicity function $n_f$ (see \S\ref{multsec}) is the absolute value of the degree of $f$ in the sense that $n_f|_I=|\deg(f|_I)|$. Thus, $n_f$ completely encodes the pm (or metric) structure of $f$.

\subsubsection{Skeletons of morphisms}
Let $f\:Y\to X$ be a generically \'etale morphism of nice curves. By a {\em skeleton} of $f$ we mean a pair $\Gamma=(\Gamma_Y,\Gamma_X)$ of skeletons of $Y$ and $X$ such that $\Gamma_Y$ contains the ramification locus $\Ram(f)$ and  $f^{-1}(\Gamma_X)=\Gamma_Y$ (in particular, $f^{-1}(\Gamma^0_X)=\Gamma^0_Y$).

\subsubsection{Semistable reduction}
It is easy to see that if $\Gamma\subseteq\Gamma'$ are two subgraphs and $\Gamma$ is a skeleton then $\Gamma'$ is a skeleton. Using this and the semistable reduction for curves one easily obtains the simultaneous semistable reduction theorem that can also be called semistable reduction of morphisms.

\begin{theor}
Any generically \'etale morphism between nice curves possesses a skeleton $\Gamma$.
\end{theor}

\begin{rem}
On the complement of a skeleton a morphism reduces to finite \'etale coverings of open discs by open discs. In general, such a morphism may have a complicated structure and this is the reason why a skeleton provides a pretty loose control on the morphism.
\end{rem}

\subsubsection{Tame morphisms}
A morphism $f$ between curves is called {\em tame} if $n_f$ takes values invertible in $\tilk$ and $f$ is called {\em wild} otherwise. A tame \'etale covering of a disc by a disc is trivial and a tame \'etale covering of an annulus by an annulus is isomorphic to the standard Kummer covering of the form $t\mapsto t^n$. So, tame morphisms are controlled by skeletons very tightly.

\subsubsection{Maps of skeletons}
More generally, it is easy to see that an \'etale covering of an annulus by an annulus is of the form $t\mapsto \sum_{i=-\infty}^\infty a_it^i$ where the series has a single dominant term $a_dt^d$ and $d>0$. In particular, the map is of degree $d$ on the skeleton. Thus, if $\Gamma$ is a skeleton of $f\:Y\to X$ then the map of graphs $\Gamma_Y\to\Gamma_X$ is enhanced to a map of metric genus graphs: to each vertex $v\in\Gamma_Y^0$ one associates a multiplicity $n_v$ and to each edge $e\in\Gamma_Y$ one associates the multiplicity $n_e$ such that $r(f(e))=r(e)^{n_e}$. These multiplicities satisfy the natural balancing conditions: if $f$ is finite then $\sum_{v\in f^{-1}(u)}n_v=\deg(f)$ for any vertex $u\in\Gamma_X^0$ and $n_v=\sum_{e\in f^{-1}(h)\cap\Br(v)}n_e$ for any vertex $v\in\Gamma_Y^0$ and an edge $h\in\Br(f(v))$ of $\Gamma_X$, where $\Br(v)$ denotes the set of all edges (or branches) coming out of $v$.

\subsubsection{Local Riemann-Hurwitz}
For a finite tame $f$ one also has the local Riemann-Hurwitz formulas: for any $v\in\Gamma_Y^0$ with $u=f(v)$ one has that $$2g(v)-2-2n_v(g(u)-1)=\sum_{e\in\Br(v)}(n_e-1),$$ which is proved by applying the RH formula to $\wt{\calH(v)}/\wt{\calH(u)}$. These formulas and the global genus formula imply the global RH formula when $X$ is proper.

\begin{rem}
One would like to extend the above formula to the non-tame case, and it is natural to expect that the local term at $e$ should be equal to the local term at the point corresponding to $e$  in the classical RH formula (e.g., see \S\ref{RHsec} below) of $\wt{\calH(v)}/\wt{\calH(f(v))}$. For non-tame morphisms two things should be modified, and we will later see that both are dealt with using the different.

(1) If $f$ is not wild at $v$ (i.e. $n_v\in\tilk^\times$) but the ramification is wild along an edge $e$ going out of $v$ then the local term $R_e$ at $e$ should be larger than $n_e-1$. So, one should naturally interpret $R_e$ in terms of the map of $k$-analytic curves.

(2) If $f$ is wild at $v$ then it often happens that $\wt{\calH(v)}/\wt{\calH(f(v))}$ is inseparable. In this case, there exists no RH-like formula based on the residue fields, and a new source of information is needed.
\end{rem}

\begin{rem}
(i) A tame $f$ is split outside of a skeleton and the only restrictions on the multiplicity function along the skeleton are the balancing conditions and the local RH formulas.

(ii) For wild maps the sets $N_{f,\ge d}$ are often huge. For example, for the Kummer map $t\mapsto t^p$ from $\bfP^1_{\bfC_p}$ to itself the set $N_{f,\ge p}$ is the metric neighborhood of $[0,\infty]$ of radius $|p|^{1/(p-1)}$.
\end{rem}

\section{The different function}
This section describes the results of \cite{CTT}. We start with recalling the definition of different and then list our main results on the different function.

\subsection{Different of extensions}

\subsubsection{The definition}
Let $L/K$ be a finite extension of real-valued fields and assume that either $K$ is discretely valued with perfect residue field or $K$ is of the form $\calH(x)$, where $x$ is a point of a nice curve. With the convention that the absolute value of an ideal $I\subseteq\Lcirc$ is $\sup_{c\in I}|c|$, the different of $L/K$ is defined to be $$\delta_{L/K}=|\Ann(\Omega_{\Lcirc/\Kcirc})|$$ if $L/K$ is separable and $\delta_{L/K}=0$ otherwise.

\begin{rem}
(i) The different measures wildness of extensions and it is multiplicative in towers.

(ii) In the case of discrete valuations one often considers the additive analogue, which is the length of $\Omega_{\Lcirc/\Kcirc}$.

(iii) In general, the different is defined using the zeroth Fitting ideal rather than the annihilator. In our case, the torsion module $\Omega_{\Lcirc/\Kcirc}$ is a subquotient of $\Lcirc$ so both definitions agree.
\end{rem}

\subsubsection{The log different}
The log different $\delta^\rmlog_{L/K}$ is defined similarly to $\delta_{L/K}$ but using the module $\Omega^\rmlog_{\Lcirc/\Kcirc}$ of logarithmic differentials. If $K$ is discretely valued then $\delta^\rmlog_{L/K}=\delta_{L/K}|\pi_L|/|\pi_K|$, and $\delta^\rmlog_{L/K}=\delta_{L/K}$ otherwise.

\subsubsection{The RH formula}\label{RHsec}
If $h\:Y\to X$ is a finite separable morphism of smooth proper connected $\tilk$-curves then the classical RH formula is
$$2g(Y)-2-2n(g(X)-1)=\sum_{y\in Y}\delta_{y/x}=\sum_{y\in Y}(\delta^\rmlog_{y/x}+n_y-1)$$
where $n=\deg(h)$, $x=h(y)$ and $\delta_{y/x}$ is the (additive) different of $k((y))/k((x))$ for $k((x))=\Frac(\hatcalO_{X,x})$ and $k((y))=\Frac(\hatcalO_{Y,y})$.

\subsection{The different function}
Let now $f\:Y\to X$ be a finite generically \'etale morphism of nice $k$-analytic curves and let $\delta_f$ be the different function introduced in \S\ref{diffunsec}.

\subsubsection{Slopes}
Naturally, $\delta_f$ contains information about the classical different at points of $Y$ of type 1 and branches of $Y$ at points of type 2. It is retrieved from the slopes (or degrees) of $\delta_f$.

\begin{theor}[{\cite[4.1.8, 4.6.4]{CTT}}]
(i) The different function extends uniquely to a pm function $\delta_f\:Y\to[0,1]$.

(ii) The slope of $\delta_f$ at a type 1 point $y$ equals $\delta^\rmlog_{y/x}$. In particular, it is positive if and only if $f$ is wildly ramified at $y$.

(iii) If $f$ is tame at a type 2 point $y$ and $v$ is a branch of $Y$ at $y$ then $\slope_v(\delta_f)=\delta^\rmlog_{v/f(v)}$.
\end{theor}

\begin{rem}
This indicates that $\delta_f$ is, in fact, the log different function. This does not affect its values at the points of $Y$ but gives a better interpretation of formulas involving differents of discretely valued fields.
\end{rem}

\subsubsection{The balancing condition}
Slopes of $\delta_f$ at a type 2 point satisfy the balancing condition of RH type which applies to the case when $\wt{\calH(y)}/\wt{\calH(f(y))}$ is arbitrary.

\begin{theor}[{\cite[4.5.4]{CTT}}]\label{balanceth}
If $y\in Y$ is of type 2 and $x=f(y)$ then $$2g(y)-2-2n_y(g(x)-1)=\sum_{v\in{\rm Br}(y)}(-\slope_v\delta_f+n_v-1).$$ In particular, almost all slopes of $\delta_f$ at $y$ equal to $n_y^i-1$, where $n_y^i$ is the inseparability degree of $\wt{\calH(y)}/\wt{\calH(f(y))}$.
\end{theor}

\begin{rem}
The balancing condition \ref{balanceth}, the formula for slopes at type 1 points, and the global genus formula imply the global RH formula when $Y$ is proper, and this can also be extended to the case with boundary. This indicates that the balancing formula is the "right" generalization of the local RH formula to the wild case.
\end{rem}

\subsubsection{The method}
The different function is a family of differents, so it is not surprising that one can describe it using a sheafified version of the definition of $\delta_{L/K}$. Namely, one considers the sheaf $\Omega^\di_X=\calO^\circ_Xd(\calO_X^\circ)$ which can be informally thought of as a version of $\Omega_{\calO_X^\circ/\kcirc}$. Then $\Omega^\di_Y/f^*\Omega^\di_X$ is a torsion sheaf of $\kcirc$-modules whose stalk at $y$ is cyclic with the absolute value of the annihilator $\delta_f(y)$. Choose $a\in\kcirc$ with $|a|=\delta_f(y)$. Reductions of $\Omega^\diamond_Y$ and $a^{-1}f^*\Omega^\diamond_X$ at $y$ induce a non-zero meromorphic map $\lam\:\tilf^*\Omega_\tilX\to\Omega_\tilY$, where $\tilf\:\tilY\to\tilX$ is the map of $\tilk$-curves associated with $\wHy/\wHx$. Then the balancing condition boils down to computing the degree of $\Omega_\tilY\otimes\tilf^*\Omega_\tilX^{-1}$ using poles and zeros of the section induced by $\lam$.

\subsubsection{The different function and the skeletons}
Let $\Gamma_Y\to\Gamma_X$ be a skeleton of $f$. It is natural to encode the balancing condition in the combinatorics of $\Gamma$. For this we should first enhance its structure by including the pm different function $\delta_\Gamma=\delta_f|_{\Gamma_Y}$. In addition, one should check whether for a vertex $y\in\Gamma_Y^0$ the skeleton contains all branches $v$ at $y$ which are {\em $\delta_f$-non-trivial}, i.e. satisfy the condition $\slope_v\delta_{f}\neq n_v-1$. It turns out that in this way one obtains a non-trivial characterization of skeletons.

\begin{theor}[{\cite[6.3.4]{CTT}}]\label{trivth}
Let $\Gamma_X$ be a skeleton of $X$ and $\Gamma_Y=f^{-1}(\Gamma_X)$. Then $(\Gamma_Y,\Gamma_X)$ is a skeleton of $f$ if and only if $\Ram(f)\subseteq\Gamma_Y^0$ and for any point $y\in\Gamma_Y$ all branches at $y$ pointing outside of $\Gamma_Y$ are $\delta_f$-trivial.
\end{theor}

\begin{rem}
(i) The behavior of $\delta_f$ completely describes the locus $N_{f,p}$ when $\deg(f)=p$. For example, if $f$ maps $\bfP^1_{\bfC_p}$ to itself by $t\mapsto t^p$ then the different is minimal and equal to $|p|$ on $[0,\infty]$ and it is trivial outside, i.e. it growths in all directions outside of the skeleton $[0,\infty]$ with slope $p-1$. This explains why $N_{f,p}$ is a metric neighborhood of radius $|p|^{1/(p-1)}$.

(ii) Even when the different $\delta_f$ behaves trivially on an interval $I=[y,z]$ its slopes depend on the multiplicity function. For example, if $n_y=p$ then the value of $\delta_f(y)$ determines $\delta_f|_I$, but if $n_y=p^n$ then to determine $\delta_f|_I$ one should know the points $x_{p^n},x_{p^{n-1}}\..x_{p^2}$ where $n_f$ drops and these points can be pretty arbitrary. In particular, the skeleton of $f$ does not control the sets $N_{f,\ge d}$ in any reasonable sense.
\end{rem}

\section{Radialization and the profile function}

\subsection{Radialization of the sets $N_{f,\ge d}$}
Let $\Gamma$ be a skeleton of a nice curve $X$, $q_\Gamma\:X\to\Gamma$ the retraction, and $r_\Gamma\:X\to[0,1]$ the inverse exponential distance from $\Gamma$. A closed subset $S\subseteq X$ is called {\em $\Gamma$-radial} if there exists a function $r\:\Gamma\to\bfR_{\ge 0}$ such that $S$ consists of all points $x\in X$ satisfying $r_\Gamma(x)\ge r(q_\Gamma(x))$.

\begin{rem}
It is easy to see that if a skeleton radializes $S$ then any larger skeleton does so.
\end{rem}

\begin{theor}[{\cite[3.3.7 and 3.3.9]{radialization}}]\label{radialth}
If $f\:Y\to X$ is a finite morphism between nice curves then there exists a skeleton of $Y$ that radializes the sets $N_{f,\ge d}$. Moreover, if $(\Gamma_Y,\Gamma_X)$ is an arbitrary skeleton of $f$ then $\Gamma_Y$ radializes these sets in each of the following cases: (1) $f$ is a normal covering (e.g. Galois), (2) $f$ is tame, (3) $f$ is of degree $p$.
\end{theor}

\begin{exam}\label{radexam}
It follows from Theorem~\ref{trivth} that if $f$ is of degree $p$ then $N_{f,p}$ is $\Gamma$-radial of radius $\delta_f^{1/(p-1)}|_\Gamma$ for any skeleton $(\Gamma,\Gamma_X)$ of $f$.
\end{exam}

\subsubsection{The splitting method}
Many results about extensions of valued fields are proved by the following splitting method:

1) Prove the result for tame extensions and wild extensions of degree $p$. Often these cases are simpler and can be managed by hands.

2) Extend the result to compositions, obtaining the case of Galois extensions.

3) Use some form of descent to deduce the non-normal case.

The splitting method extends to a local-analytic setting because the category of \'etale covers of a germ $(X,x)$ of an analytic space at a point is equivalent to the category of \'etale $\calH(x)$-algebras by a theorem of Berkovich. Theorem~\ref{radialth} is proved easily by the splitting method since the tame case is clear and the degree-$p$ case is controlled by the different by Example~\ref{radexam}.

\subsection{The profile function}
One may wonder if the radii of the sets $N_{f,p^n}$ are reasonable functions analogous to the different. The answer is yes, but the best way to work with them is to combine them into a pm function from $[0,1]$ to itself.

\subsubsection{$\Gamma$-radial morphisms}
Let $f\:Y\to X$ be a morphism and $\Gamma=(\Gamma_Y,\Gamma_X)$ a skeleton of $f$. For a point $a\in Y$ of type 1 consider the interval $I=[a,q_{\Gamma_Y}(a)]$ and identify it with $[0,1]$ via the radius parametrization. Similarly, identify $f(I)=[f(a),q_{\Gamma_X}(f(a)]$ with $[0,1]$. Then $f|_I$ is interpreted as an element of $P_{[0,1]}$ and we say that $f$ is {\em $\Gamma$-radial} if $f|_I=\phi_q$ depends only on $q=q_{\Gamma_Y}(a)$. In this case we say that $\phi\:\Gamma_Y\to P_{[0,1]}$ is the profile function of $f$.

\begin{rem}
(i) It is easy to see that $f$ is $\Gamma$-radial if and only if all sets $N_{f,d}$ are $\Gamma_Y$-radial and then the breaks of $\phi_q$ occur at the radii of the sets $N_{f,p^n}$.

(ii) Thus, the radialization theorem implies that any finite morphism is $\Gamma$-radial for a large enough skeleton $\Gamma$. In particular, this gives another way to define $\phi_q$: it is the map $f|_I$ for a generic interval connecting $q$ to a type 1 point.

(iii) The profile function is a much more convenient invariant than the set of radii of $N_{f,d}$, mainly because it is compatible with compositions of radial morphisms.
\end{rem}

\subsubsection{Interpretation as Herbrand function}
It turns out that the profile function can be interpreted using a classical invariant from the theory of valued fields. It is well-known that for a finite seprable extension $l/k$ of discrete valuation fields with perfect residue fields, the Herband function $\phi_{l/k}$ is a multiplicative (with respect to towers of extensions) invariant which efficiently encodes nearly all information about the wild ramification properties of $l/k$.

It is shown in \cite[\S4]{radialization} that the theory extends to extensions of the form $\calH(y)/\calH(x)$, where $y,x$ are points on $k$-analytic curves (for an algebraically closed $k$). The only technical obstacle is that in the classical theory one crucially uses that the extension of integers is monogeneous while $\calH(y)^\circ/\calH(x)^\circ$ is integral but does not have to be finite. However, it is shown in \cite[4.2.8]{radialization} that $\calH(y)^\circ/\calH(x)^\circ$ is almost monogeneous in the sense that $\calH(y)^\circ$ is a filtered union of subrings of the form $\calH(x)^\circ[t]$, and it is shown in \cite[\S4.3]{radialization} that the classical theory extends to extensions of analytic fields with almost monogeneous extensions of rings of integers.

\begin{theor}[{\cite[4.5.2]{radialization}}]\label{herbrandth}
If $f\:Y\to X$ is a generically \'etale morphism between nice curves then for any point $y\in Y$ of type 2 with $x=f(y)$ the profile function $\phi_y$ coincides with the Herbrand function $\phi_{\calH(y)/\calH(x)}$ of the extension $\calH(y)/\calH(x)$.
\end{theor}

\begin{rem}
(i) The proof is via the splitting method using that for extensions of degree $p$ the Herbrand function is determined by the different (it has slopes 1 and $p$ and the break point is determined by the different).

(ii) The theorem gives a natural geometric interpretation of Herbrand function which works for all extensions (even inseparable ones) on the equal footing. Note that the classical Herbrand function is defined first for Galois extensions and then extended to arbitrary ones by multiplicativity.
\end{rem}

\subsubsection{Piecewise monomiality of the profile function}
It is natural to expect that the profile functions $\phi_y$ should discover a nice global behavior. Indeed, one can easily introduce a notion of pm functions on $Y$ with values in $P_{[0,1]}$ and the following result holds.

\begin{theor}[{\cite[3.4.8]{radialization}}]
If $f$ is as above then the family of profile functions $\phi_y$ extends uniquely to a pm function $\phi\:Y^\hyp\to P_{[0,1]}$.
\end{theor}

\bibliographystyle{amsalpha}
\bibliography{wild_covers}

\end{document}